\documentclass [10pt,twoside]{article} 
\usepackage{times} 
\pdfoutput=1

\usepackage{amsmath}
\usepackage{amsthm}
\usepackage{amsfonts}
\usepackage{txfonts}
\usepackage{amssymb}
\usepackage{amsmath}
\usepackage{latexsym}
\usepackage{multirow}
\usepackage{slashbox}

\usepackage[all]{xy}
\usepackage{epsfig}   
\usepackage{color}
\usepackage{ifthen}
\usepackage{verbatim}
\usepackage[numbers, sort&compress]{natbib} 

\def\bR{\begin{color}{red}} 
\def\bB{\begin{color}{blue}}
\def\bM{\begin{color}{magenta}}
\def\bC{\begin{color}{cyan}}
\def\bW{\begin{color}{white}}
\def\bBl{\begin{color}{black}} 
\def\bG{\begin{color}{green}}
\def\bY{\begin{color}{yellow}}
\def\e{\end{color}}

\theoremstyle{definition}
\theoremstyle{definition}

\newcommand{\bit}{\begin{itemize}}
\newcommand{\eit}{\end{itemize}\par\noindent}
\newcommand{\ben}{\begin{enumerate}}
\newcommand{\een}{\end{enumerate}\par\noindent}
\newcommand{\beq}{\begin{equation}}
\newcommand{\eeq}{\end{equation}\par\noindent}
\newcommand{\beqa}{\begin{eqnarray*}}
\newcommand{\eeqa}{\end{eqnarray*}\par\noindent}
\newcommand{\beqn}{\begin{eqnarray}}
\newcommand{\eeqn}{\end{eqnarray}\par\noindent}

\def\rto{\multimap}
\def\lto{\multimapinv}

\theoremstyle{plain}
\theoremstyle{plain}
\newtheorem{thm}{Theorem}
  \newtheorem{defn}[thm]{Definition}

\newcommand{\xyR}[1]{%
\xydef@\xymatrixrowsep@{#1}}
\newcommand{\xyC}[1]{%
\xydef@\xymatrixcolsep@{#1}}

\title{An alternative Gospel of  structure:\\ order, composition, processes}
\author{Bob Coecke\\ 
University of Oxford,\\ Department of Computer Science\\
coecke@cs.ox.ac.uk}
\date{}

\begin{document}
\maketitle
\thispagestyle{empty}

\begin{abstract}
We survey  some basic mathematical structures, which  arguably are more primitive than the  structures taught at  school.   These structures are orders, with or without composition, and  (symmetric) monoidal categories. We list several `real life' incarnations of each of these.  This chapter also serves as an introduction to these structures and their current and potentially future uses in linguistics,  physics and knowledge representation.  

This paper is the introductory chapter of the volume \cite{FlowinCats}, which collects a range of tutorial papers spanning a wide range of subjects employing structures of the kind discussed here, authored by leaders of their respective fields.
\end{abstract}


\bigskip 
What are the  fundamental  mathematical structures? Evidently, in order to address this question  one first needs to address another question, namely,  the question on  what a mathematical structure actually is.  There are several options available.

Are mathematical structures the objects of mathematical theories?  And is it then  the mathematician  who decides which ones are truly fundamental? It is indeed often said by mathematicians that good mathematics is what good mathematicians do.   It has also been strongly argued by mathematicians,  that good mathematics should be a discipline which exists in isolation from the (other) natural sciences, e.g.~by Bourbaki, and that the use of example applications and illustrations outside of straight mathematics is to be discouraged.  We find these views somewhat circular and solipsistic, and even disrespectful to the other disciplines,  in the light of the origins and history of mathematics. 

From an alternative more reductionist perspective, one may think that  the  fundamental  mathematical structures are the simple things from which one can build more complex things.   Proposed candidates of fundamental structures include the  Platonic solids from which the \em classical elements \em (Earth, Water, Air, and Fire) were constructed, and more recently, \em sets\em, which within the realm of set-theory  are supposed to underpin all of mathematics. In our digital age, 0's and 1's underpin all of our data.  So do 0's and 1's constitute the appropriate language of communicating our opinions about the movies that are encoded in terms of 0's and 1's on a computer's hard drive?  Of course not.  The fact that one can `code' everything by means of one particular structure, just as most modern mathematics can be encoded in sets, does not mean that it is the most convenient language for discourse. 
 

Or, are fundamental  mathematical structures  the objects which underpin the world in which we live, i.e.~our reality? This raises the question, what do we mean by reality? Is this restricted to material reality or does this also include the stuff going on in our minds?  In Aristotelianism numbers exist as concepts, which means that they do exist, but only in our mind. So  where were those numbers  before there were any minds that they could inhabit.  
Today, number-based calculus is still a pillar of modern mathematics, and many other fields of mathematics have some sort of numbers (be it either integer, real or complex) as their starting point.  
At school, first we learn natural numbers, then negatives, then rationals, reals, and eventually complex numbers.  We can organize these as matrices yielding vector spaces, which is probably the most prominent mathematical structure at the moment. This points at these number systems maybe being the fundamental mathematical structures.

In this paper we will preach an alternative Gospel, which qualifies structures in terms of their high-levelness in use.    As a starting point, we try to analyze several forms of practice, starting in fact with `counting' as underpinning trade, as well as basic mathematical practice and its relation to speech and reasoning, as well as nature.  This is by no means a thorough development but just a few stone throws towards trying to develop a story of this kind.

Our preliminary development  involves  three stages:
\bit
\item[(1)] Firstly, we considered ordering of things.
\item[(2)] Secondly, we assume that these things can be composed.
\bit
\item[(2')] These ordered composable things also provide a stepping stone to the next stage. We can consider processes on a single thing that can be composed in time, and that moreover happen to be ordered.  Dropping this ordering yields a stage (2') which is complementary to (2).  Together, stages (2) and (2') combine into the remaining stage, that is:
\eit
\item[(3)] We consider processes, composable in time,  for composable things.  
\[
\xymatrix@=0.20in{
& (1) \ar[ld]\ar[rrd] && \\
(2)\ar[rd] && (2')\ar[ld]\ar[r] & (2'+1)\\
& (3)
}
\]
\eit
We will make a case for why we consider these structures as fundamental.

For each of these stages we provide a number of example applications in different areas, which indicate the general discipline transcending nature of theses structures.  We consider this as additional evidence of the foundational nature of these structures.

But as a matter of convenience, we will still rely on set theory to provide us with a notion of equality in the formal definitions, as this will appeal to what the reader is used too. And of course, most existing example applications that we want to present have been formulated in set-theoretic terms, so we don't really have much choice here.  However, it should be clear that these examples  have a `meaning in reality' which does not require this particular underpinning.  We will illustrate this for one specific example. Consider the case of describing the possible states of three things `in set theory vs.~in reality'.  

In set theory, each of these things will be some structured set, that is, a set $X$ together with additional operations which may encode topology (a collection of open subsets), geometry (lines etc.), or algebraic structure (e.g.~a binary operation).   The operation `putting  things together' would be encoded in most cases as a cartesian product.  Focussing from now on only on the supporting sets, we can build a triple of things in three manners:
\bit
\item First we combine $X$ and $Y$ into $X\times Y:= \{(x, y)\mid x\in X, y\in Y\}$ and then we combine $X\times Y$ and $Z$ into $(X\times Y)\times Z:= \{((x, y), z)\mid x\in X, y\in Y, z\in Z\}$.
\item First we combine $Y$ and $Z$ into $Y\times Z:= \{(y, z)\mid y\in Y, z\in Z\}$ and then we combine $X$ and $Y\times Z$ into $X\times (Y\times Z):= \{(x, (y, z))\mid x\in X, y\in Y, z\in Z\}$.
\item We do everything at once by considering triples $\{(x, y, z)\mid x\in X, y\in Y, z\in Z\}$.
\eit
Now, the three results are not equal!  Instead, they are isomorphic in a very `natural' manner, since we can pass from $\{((x, y), z)\mid x\in X, y\in Y, z\in Z\}$ to $\{(x, (y, z))\mid x\in X, y\in Y, z\in Z\}$ by `re-bracketting', and to  $\{(x, y, z)\mid x\in X, y\in Y, z\in Z\}$ by `dumping brackets'.  

Clearly `re-bracketting' and  `dumping brackets' are operations that have many nice properties, which boil down to saying that each of these three ways of combining three things are essentially the same for all practical purposes.   But to properly state this in mathematical terms requires some pretty heavy machinery.  For example, simply saying that these are isomorphic is not enough, since along such an isomorphism we may loose the fact that we truly are considering three things, e.g.
\[
\{0,1\}\times\{0,1\}\times\{0,1\}\simeq\{0,1,2,3,4,5,6,7\}\ .
\]

To capture the `natural' manner in which constructions like
\[
\{((x, y), z)\mid x\in X, y\in Y, z\in Z\}\ \  \mbox{\rm and}\ \  \{(x, (y, z))\mid x\in X, y\in Y, z\in Z\}
\]
are equivalent Eilenberg and MacLane introduced category theory, which first required to define categories, then functors, and finally natural transformations \cite{EilenbergMacLane}, a highly non-trivial concept.  For something as simple as saying that there are three things,  these concepts feel a bit like overkill.  Indeed, in reality, given three objects, it simply doesn't matter that we first consider two of them together, and then the third one, or the other way around, in order to describe the end result.  In other words, the bracketing in the above description has no counterpart in reality, it is merely imposed on us by set theory!  

But make no mistake, it is not category theory that we blame here in any manner, but set theory.  Certain branches of category theory have no problem in describing three things on an equal footing, in terms of so-called strict symmetric monoidal categories.  True, in its usual presentations (e.g.~\cite{MacLane}) it is build upon set theory, but not to do so is exactly what we are arguing for here.   An elaborate discussion of this point can be found in \cite{CatsII}. Also, category theory is a huge area and by no means all of it will be relevant here.  

\section{Order}  

Above we mentioned that number-based calculus is a pillar of modern mathematics.

But probably more primal than counting is the simple fact wether one number is \em larger than \em another.  If numbers characterize the cost of something, say a piece of meat, being able to afford that piece of meat requires one to have  \em more \em money available \em than \em the cost of that piece of meat.  The meaningfulness of the question ``Can I afford this piece of meat?'' results from the fact that numbers are \em ordered\em: if $x$ is the price and $y$ is one's budget, then `Yes I can!'' means $x\leq y$ while ``No I cannot!'' means $x\not\leq y$.   

\begin{defn}\em 
A \em total ordering \em on a set $X$ is a relation $\leq$ on $X$, i.e.~a collection of pairs $(x, y)\subseteq X\times X$ which is \em anti-symmetric\em, \em transitive \em and \em total\em, that is, respectively,
\bit
\item  $\forall x, y\in X: x\leq y, y\leq x \Rightarrow  x= y$
\item $\forall x, y, z\in X: x\leq y, y\leq z \Rightarrow  x\leq z$
\item $\forall x, y, z\in X: x\leq y$ or $y\leq x$
\eit
\end{defn}

Of course, one cannot compare apples and lemons, if we don't consider monetary   value.  Three apples are less than four apples, but what about three apples versus three lemons?  

\begin{defn}\em\label{def:poset}
A \em preordering \em  on $X$ is a transitive relation $\lesssim$ which is also \em reflexive\em, that is,   
\bit
\item  $\forall x\in X: x\lesssim x$
\eit
A \em partial ordering \em is an anti-symmetric preordering, and denoted $\leq$. 
\end{defn} 

By definition, each partial ordering is a preordering,  and also, each total ordering is a partial ordering since totality implies reflexivity (apply totality when $x=y$).  

Moreover, each preordering yields an \em equivalence relation\em, which is defined as   a preordering $\simeq$ which is  also \em symmetric\em, that is, 
\bit
\item  $\forall x, y\in X: x\simeq y \Leftrightarrow y\simeq x$\,,
\eit
simply by setting $x\simeq y$ to be  $x\lesssim y$ and $y\lesssim x$.  Then, the corresponding set of equivalence classes $\{C_x\mid x\in X\}$, where $C_x=\{y\in X\mid y\simeq x\}$,   forms a partial ordering when setting $ C_x\leq C_y$ whenever $x\lesssim y$.  For example, if we order things in terms of their cost we obtain a preordering, since there may be many things that have the same cost.  This gives rise to a partial ordering, which is in fact a total ordering, of the occurring costs themselves.

Sometimes one encounters the notion of a \em strict partial order \em $<$, which is \em irreflexive \em ($x\not\leq x$), \em asymmetric \em ($a< b\Rightarrow b\not <a$) and transitive, but it is easily seen that each partial ordering as in Defn.~\ref{def:poset} induces a strict partial ordering by setting $x<y$ when $x\leq y$ and $x\not =y$; conversely, $x\leq y$ when  $x<y$ or $x=y$ turns a strict partial ordering into a partial ordering as in Defn.~\ref{def:poset}. So both notions are equivalent in their uses.

Preorderings and partial orderings arise in many everyday situations, as well as in important areas of science.  Arguably, they are the most primitive mathematical concept imaginable, capturing the very basic notion of `comparing'.   

\subsection{Reachability, causality and relativity}

Consider a collection of events, a certain rock  concert in New Orleans on a particular day, the marriage  of some old friend, Carnival in Rio this year, etc. Then there exists a partial ordering where the relation $\leq$ stands for: ``if I attend event $a$, can I also attend event $b$?''  Otherwise put, the ordering $a\leq b$ captures wether one can reach event $b$ from event $a$.  The overall data in this partial ordering  is imposed by the transport network of the world.  

In physics, due to the velocity bound imposed by the speed of light, a similar partial order exists which encodes wether light can travel from one point $x$ in space-time to another point $y$ in space-time.  In the case of yes we can write $x\leq y$. And since nothing travels faster than time, this partial order encodes which event in space-time can causally affect which other event in space-time.  In fact, there exist results that enable one to reconstruct the entire geometric space-time manifold from this partial ordering e.g.~\cite{MartinPanagaden}.  There are moreover several research programs that take partial ordering not only as a framework for discussing special and general relativity, but also as the basis for crafting a theory of quantum gravity \cite{Bombelli, Sorkin}, the Holy Grail of modern physics.

Also in computer science, identical ideas on causality as partial orderings exist \cite{Lamport}, and through a pair of  `partial order'-glasses, both the areas of relativity theory and the organization of events in a distributed computational system look remarkably similar. In other words, in high-level terms they essentially coincide.

\subsection{Information content,  propositional knowledge and domains}

One may be interested in the information content of some pieces of data, for example in terms of entropy. Typical measures of information content will assign some real number, so that, since real numbers are totally ordered, we can decide which one of two of pieces of data is the most informative one.  This then induces a total preordering on these pieces of data.  In many cases, the actual values don't really matter, merely what is more informative than what matters.  So in fact, rather than taking a valuation into real numbers, one can just consider a total preordering on the pieces of data.

Now of course, if two pieces of data each constitute one bit, while their information content is the same, these may be incomparable in terms of their propositional content. For example, while ``$x$ is an apple'' implies ``$x$ is fruit'', it is incomparable with ``$x$ is a lemon''.  As propositions, these pieces of data form a partial order, where $a\leq b$ stands for the fact that $a$ implies $b$, e.g.~``being an apple'' implies ``being fruit''.  This idea is the cornerstone to \em algebraic logic \em \cite{Davey}, discussed in the next section.  

For the purpose of combining informative and propositional content, domain theory was crafted \cite{ScottDomains}, as a new mathematical foundation for computation.  Here, domains are partial orders in which certain subsets have least upper bounds.

\begin{defn}\em
Given a partial ordering $(X, \leq)$, a subset $Y\subseteq X$ has a   \em least upper bound \em (or \em join\em) if there is an element $x \in X$ which is such that:
\bit
\item $\forall y\in Y: y\leq x$
\item $\forall x'\in X: (\forall y\in Y: y\leq x')\Rightarrow x\leq x'$\ .
\eit
In this case we denote this element $x$ as $\bigvee Y$. A \em greatest lower bound \em (or \em meet\em)  $\bigwedge Y$ of $Y\subseteq X$ is defined similarly simply by replacing $\alpha \leq \beta$ by $\beta \leq \alpha$ in the above.
\end{defn}

While we won't discuss the nature of those subsets that have least upper bounds in a domain, a particular example of a least upper in algebraic logic  is \em disjunction \em $a\vee b$, where $Y$ consists of two elements. \em Conjunction  \em $a\wedge b$ is the corresponding greatest lower bound.

Interestingly, the partial orders underpinning space-time are in fact also domains \cite{MartinPanagaden}. Hence at this order-theoretic level again two seemingly disjoint subjects become  the same when taking a sufficiently high-level perspective.  

Also in this context, while it fails to be a domain and only captures propositional content up to symmetries, the \em majorization preordering \em \cite{Muirhead} on probabilities has a range of applications, ranging from economy \cite{MajorizationBook} to quantum information theory  \cite{Nielsen}, where it captures degrees of entanglement.  What is ordered here are descending discrete probability distributions, that is, $n$-tuples $(x_1, \ldots, x_n)$ with $\sum_i x_i = 1$ and that are such that $x_i\geq x_{i+1}$. We say that an ordered $n$-tuple $(x_1, \ldots, x_n)$ is \em majorized \em by  $(y_1, \ldots, y_n)$ if 
\[
\forall k\in \{1, \ldots, n-1\}: \sum_{i=1}^{i=k} x_i \leq \sum_{i=1}^{i=k} y_i\ .  
\]
Intuitively, this means that $(y_1, \ldots, y_n)$ is a `narrower' or `sharper' distribution than $(y_1, \ldots, y_n)$.

Unfortunately, majorization does not extend to a partial order on all probabilities, so it fails to capture propositional content. A genuine partial ordering on probabilities in which propositional structure naturally embeds   is in \cite{CoeckeMartin}.  It is `desperately looking for more applications', so please let us know if you know about one!

We refer the reader to \cite{AbramskyJung} for the role of domain theory in computer science, where it originated,  and in particular to Martin's tutorial  \cite{MartinTutorial} for a much broader range of applications in a variety of disciplines.

\subsection{Logic and the theory of mathematical proofs}\label{sec:logicandproof}

In algebraic logic, which traces back to Leibniz,  one typically would like to treat ``$a$ implies $b$'' itself as a proposition, something which is realized by an \em implication connective\em, that is, an operation $\Rightarrow: X\times X \to X$ on the partial ordering $(X,\leq)$.  Typically, one also of course assumes conjunction and disjunction, and in its weakest form implication can be defined by the following order-theoretic stipulation:
\beq\label{eq:Heyting}
(a \wedge b) \leq c\ \  \mbox{if and only if}\ \  a\leq (b\Rightarrow c)
\eeq
One refers to partial orderings with such an implication as \em Heyting algebras\em.  From eq.~(\ref{eq:Heyting}) it immediately follows that  the \em distributive law \em holds, that is, 
\[
a \wedge (b \vee c)= (a \wedge b) \vee (a \wedge c)\ , 
\]
or in terms of words, 
\[
a \ AND\ (b \ OR \ c)= (a\ AND\ b)\ OR\ (a\ AND\ c) \ .
\]

A special case of implication is a \em Boolean implication\em, which is defined in terms of negation as $a\Rightarrow b := \neg a \vee b$.  We won't go into the details of defining the negation operation $\neg: X\to X$, but one property that it has is that it reverses the ordering, that is:
\bit 
\item $a\leq b\  \Leftrightarrow \ \neg b\leq\neg a$\ .
\eit
In the case of a Boolean implication we speak of a \em Boolean algebra\em.  The smallest Boolean algebra has two elements, $true$ and $false$, with $false\leq true$.

Closely related to algebraic logic is \em algebraic proof theory\em, a meta-theory of mathematical practice.  One may indeed be interested from which assumptions one can deduce which conclusions, and this  gives rise to a preordering, where $a\lesssim b$ stands for `from $a$ we can derive $b$'.  Since possibly, from $a$ we can derive $b$, as well as, from $b$ one can derive $a$, we are dealing with a proper preorder rather than a partial order. In this context, eq.~(\ref{eq:Heyting}) corresponds to the so-called \em deduction theorem \em \cite{Kleene}.  It states that if one can derive $c$ from the assumptions $a$ and $b$, then this is equivalent to deriving that $b$ implies $c$ given $a$.  

We stress here that $\lesssim$ now represents the \em existence \em of a proof, but not a proof itself.  We will start to care about the proofs themselves in Sec.~\ref{sec:SMC}.

\subsection{General `things' and processes}\label{sec:things}

The idea of the existence of a mathematical proof that transforms certain assumptions into certain conclusions extends to general situations where processes transform certain things into other things.  For example, if I have a raw carrot and a raw potato I can transform these into  carrot-potato mash.  On the other hand, one cannot transform an apple into a lemon, nor  carrot-potato mash into an egg.  So in general, when one considers a collection of things (or systems), and processes which transform things into other things, one obtains a preordering that expresses what can be transformed into what.  More precisely, $a\leq b$ means that \em there exists a process that transforms $a$ into $b$\em.  The technical term for things in computer science would be \em data-types\em, and for processes \em programs\em, in physics things are \em physical systems \em and examples of processes are evolution and measurement, in cooking the thing are the \em ingredients \em while example processes are boiling, spicing, mashing, etc.

 \begin{center} 
\begin{tabular}{c|c|c|} 
& thing / system & process  \\ \hline \hline
Math.~practice & propositions & proofs (e.g.~lemma, theorem, etc.) \\ \hline
Physics & physical system & evolution, measurement etc. \\ \hline
Programming & data type & program  \\ \hline
Chemistry & chemical & chemical reaction \\ \hline
Cooking & ingredient & boiling, spicing, mashing, etc. \\ \hline
Finance & e.g.~currencies & money transactions \\ \hline
Engineering & building materials & construction work \\
\end{tabular}
\end{center} 

\section{Orders and composition}\label{sec: composition}

We mentioned that one cannot compare apples to lemons.  However three apples and two lemons is clearly less that five apples and four lemons.  To formalize this, we need to have a way of saying that we are `adding apples and lemons'. We will refer to this `adding' as composition.  Clearly, this composition needs to interact in a particular way with the ordering such that, either increasing the apples or the lemons increases the order of the composite. Note that adding apples to apples, or money to money, is just a special case of this, where the ordering is total, i.e.~everything compares to everything, rather than being a proper partial ordering in which some things don't compare. 

\begin{defn}\em
A \em monoid \em is a set $X$ together with a binary operation $\cdot : X\times X\to X$ which is both associative and admits a two-sided unit $1\in X$, that is, respectively:
\bit
\item $\forall x, y, z\in X:  x\cdot (y \cdot z)= (x\cdot y) \cdot z$
\item $\forall x\in X: 1\cdot x = x\cdot 1 = x$
\eit
A \em totally/partially/preordered monoid \em is a set $X$ which is both a monoid and a total/partial/pre-order, and which moreover satisfies \em monotonicity of the monoid multiplication\em:
\bit
\item $\forall x, y, x', y' \in X: x\lesssim y, x'\lesssim y' \Rightarrow x\cdot x' \lesssim y\cdot y'$
\eit
\end{defn}

Many of the applications mentioned above extend to ordered monoids. For example, in algebraic logic, either the conjunction or the disjunction operation yields a partially  ordered monoid.  When tinking about things and processes, we can obtain an ordered monoid as  in the case of apples  and lemons. We can compose different things, and we can compose the processes acting thereon in parallel.  In the case of mathematical proofs this would simply mean that we prove $b$ given $a$ as well as proving $d$ given $c$.

While in the case of apples and lemons composition is \em commutative\em, that is:
\bit
\item $\forall x, y\in X:  x\cdot y= y \cdot x$\,, 
\eit
composition of processes \em in time \em is typically non-commutative.  

And we can indeed also think of the elements of an ordered monoid as processes themselves, the monoid composition then standing for process $b$ \em after \em process $a$.  The ordering is then an ordering on processes, for example, given a certain proof, the proof of either a stronger claim from the same assumptions, or the same claim from weaker assumptions  would be strictly below the given proof. 

Summing the above up we obtain two distinct ordered monoids for the example of mathematical proofs (as well as for the other examples):
\bit
\item {\bf Composition of things}, or {\bf parallel composition}, where the ordered propositions can be composed, that is, we can consider collections of assumptions rather than individual ones, and comparing these may happen component-wise, but doesn't have to (cf.~conjunction and disjunction as composition).
\item {\bf Composition of processes}, or, {\bf sequential composition}, where the elements of the monoid are the proofs themselves, rather than propositions. The composition is now `chaining' proofs, that is, a proof of $b$ given $a$ and a proof $c$ given $b$ results in a proof of $c$ given $a$.  The ordering is then in a sense  the quality of the proof, in that a better proof is one that achieves stronger conclusions from weaker assumptions.
\eit
More generally, non-commutative ordered monoids yield some interesting new applications, for example, in \em reasoning about knowledge  \em and in \em natural language\em.

\subsection{Galois adjoints as assigning causes and consequences} 

\begin{defn}\em
For two order-preseving maps $f:A\to B$ and $g: B\to A$ between partially ordered sets $(A,\leq)$ and $(B,\leq)$ we say that $f$ is \em left adjoint \em to $g$ (or equivalently, that $g$ is \em right adjoint \em to $f$), denoted $f \dashv g $, if we have that: 
\bit
\item $\forall a\in A, b\in B: f(a)\leq b \Leftrightarrow a \leq g(b)$\ ,  
\eit
or what is easily seen to be equivalent \cite{Koslowski},  if we have that:  
\bit
\item $\forall b\in B: f(g(b))\leq b$ and $\forall  a\in A: a\leq g(f(a))$\ . 
\eit
\end{defn} 

This at first somewhat convoluted looking definition has a very clear interpretation if we think of $f$ as a process which transform propositions $a\in A$ of system $A$ into propositions $b\in B$  of $B$.  This goes as follows.   

Assume we know that process $f$ will take place (e.g.~running it on a computer as a computer program) but we would want to make sure that after running it $b\in B$ holds, and the means we have to impose this is to make sure that before running it some $a\in A$ holds. Then, the Galois adjoint to $f$ gives us the answer, namely, the necessary and sufficient condition is to take $a\leq g(b)$.  In computer science one refers to $ g(b)$ as the \em weakest precondition \em to realize $b$ by means of $f$ \cite{Dijkstra,Hoare}.  More generally, one can think of $g$ as assigning \em causes\em, while $f$ assigns \em consequences \em for any kind of process \cite{CMS}.

Note in particular  that eq.~(\ref{eq:Heyting}) is also  of the  form of a Galois adjunction.  Indeed, explicitly putting quantifiers and re-ordering symbols, eq.~(\ref{eq:Heyting}) can be rewritten as:
\bit
\item $\forall c\in X, \bigl(\forall a\in X, b\in X: (c\wedge a) \leq b \Leftrightarrow a \leq (c\Rightarrow b)\bigr)$\ .
\eit
The expression within the large brackets is  a Galois adjunction for $A=B=C:=X$ and: 
\[
f:=(c\wedge -):X\to X\quad \mbox{\rm and} \quad g:=(c\Rightarrow -):X\to X\ . 
\]
In eq.~(\ref{eq:Heyting}) we ask this  Galois adjunction to be true for all $c\in X$, and all these conditions together then define an implication connective $(- \Rightarrow -):X\times X\to X$.  More generally, adjointness provides a comprehensive foundation for logic \cite{LawvereAdjointness,LambekScott}.

While the author does not subscribe (anymore) to this approach, so-called \em quantum logic \em \cite{BvN, Piron} is also an order-theoretic enterprise, and the corresponding `weak' implication, the Sasaki hook \cite{Sasaki}, can best be understood in terms of Galois adjoints.  

Birkhoff and von Neumann noted that observationally verifiable propositions of a quantum system still form a partial ordering, which admits least upperbounds and greatest lower bounds, but that the distributive law fails to hold, and as a consequence, that there is no connective $(-\Rightarrow-)$ as in e.q.~(\ref{eq:Heyting}) since that would imply distributivity (see Sec.~\ref{sec:logicandproof}).  However, there is an operation $(c\leadsto-): X\times X\to X$ for all $c\in X$ which is such that:
\bit
\item $\forall a\in X, b\in X: P_c(a) \leq b \Leftrightarrow a \leq (c\leadsto b)$\ .
\eit
where $P_c:X\to X$ stands for the orthogonal  projection on $C$, or in physical terminology, the collapse onto $c$.  The collapse is an actual physical process that happens when measuring a quantum system, so in the light of the above discussion on causes and consequences, which as discussed in \cite{CMS} extends to these quantum logics,  the operation $(c\leadsto -): X\to X$ should be understood as assigning the weakest precondition $(c\leadsto b)\in X$ that has to hold before the collapse in order for $b\in X$ to hold after the collapse.

While as already mentioned, the author does not subscribe to quantum logic anymore, Constantin Piron's operational take on the subject \cite{Piron, DJMoore} has greatly influenced the author's thinking.  Unfortunately, Piron died during the final stages of writing of this chapter.

\subsection{Dynamic (\&) epistemic logic}

Above we briefly discussed propositional logic, that is, describing the set of propositions about a system as a partially ordered set. So what about the change of these propositions?  \em Actions\em, which change propositions can be described as maps acting on these propositions, and form themselves an ordered monoid.  It is easily seen that these maps should preserve disjunctions, in order theoretic terms that is, greatest lower bounds.  Indeed, if $a\ OR \ b$ holds, than after the action $f$, clearly $f(a)\ OR \ f(b)$ should hold, so  $f(a\vee b)=f(a)\vee f(b)$.   This guarantees that these actions have left Galois adjoints, assigning causes.

Now, while propositions may be as they are, an agent may perceive them differently.  This may  for example  be due to \em lying actions \em by some agent who is supposed to communicate these changes of propositions.  These situations were considered by Baltag, Moss and Solicki in \cite{BMS}, to what is referred to as \em dynamic epistemic logic\em, and it was shown in \cite{BCS2} that all of this is most naturally cast in order theoretic terms.  

This setting comprehends and stretches well beyond the fields of epistemic logic \cite{Wiebe} and dynamic logic \cite{Harel}, both conceptual variations on so-called \em modal logic \em \cite{Kripke}, and all of which are part of modern algebraic logic. They have a wide range of applications in computer science, including soft- and hardware verification \cite{Sterling}.  On the other hand, these logics also underpin Carnap's philosophy on semantics and ontology \cite{Carnap}.

\subsection{Linguistic types}\label{sec:catgram} 

We can also compose words in order to build sentences.  However, not all strings of words make up meaningful sentences, since meaningfulness imposes constraints on the grammatical types of the words in the sentence.  

By having a partial order relation besides a composition operation we can encode how the overall grammatical type of a string of words evolves when composing it with other types, and ultimately, make up a sentence.  The ordering $a_1\cdot\ldots\cdot  a_n \leq b$ then encodes the fact that the string of words of respective grammatical  types $a_1$, \ldots, $ a_n$   has as its  overall type $b$.   For example, if $n$ denotes  the type of a noun,  $tv$ the type of a transitive verb and $s$ the type of a (well-formed) sentence,  then $n\cdot tv \cdot n\leq s$ expresses the fact that a noun (= object), a transitive verb, and another noun (= subject), make up a well-formed sentence.  

The remaining question is then how to establish a statement like $a_1\cdot\ldots\cdot  a_n \leq b$.  The key idea is that some grammatical types are taken to be atomic, i.e.~indecomposable, while others are compound, and one considers additional operations, which may either be unary or binary, subject to some laws that allow one to reduce  type expressions.  For example, assuming that one has left- and right-`pre-inverses', respectively denoted ${}^{-1}(-)$ and $(-)^{-1}$, and subject to $a\cdot {}^{-1}a\leq 1$ and $a^{-1}\cdot a\leq 1$, then for the compound transitive verb type: 
\[
tv= {}^{-1}n\cdot s\cdot n^{-1}
\]
we have:  
\[
n\cdot tv\cdot n= n\cdot {}^{-1}n\cdot s\cdot n^{-1}\cdot n\leq 1\cdot s\cdot 1\leq s\ ,
\]
so we can indeed conclude that $n\cdot tv \cdot n$ forms a sentence.

As the area of type grammars hasn't reached a conclusion yet on the question of which partially ordered monoids, or in short, pomonoids, best captures `universal grammatical structure', we give a relatively comprehensive historical overview of the structures that have been proposed.  Moreover, several of these will provide a stepping stone  for the categorical structures in the next section.  Historically, the idea of universal grammar of course traces back to Chomsky's work in the 50's \cite{Chomsky}. The mathematical development has mainly be driven by Lambek \cite{Lambek0, Lambek1, LambekBook, MoortgatBook}, in many stages spanning some 60 years of work.
  
\begin{defn}\em
A \em protogroup \em \cite{Lambek1} is a pomonoid  
\[
(X, \leq, {}^*(-), (-)^*) 
\]
where ${}^*(-):X\to X$ and $(-)^*:X\to X$  are such that:
\bit
\item $\forall a, b\in X:   a\cdot {}^*a \leq 1$ and $b^* \cdot b\leq 1$\ .
\eit
\end{defn}

\begin{defn}\em
 An \em Ajdukiewicz-Bar-Hillel pomonoid \em \cite{Ajdukiewicz,Bar-Hillel} is a pomonoid  
\[
(X, \leq, (-\rto -), (-\lto-)) 
\] 
where $(-\rto -):X\times X\to X$ and $(-\lto -):X\times X\to X$  are   such that:
\bit
\item $\forall a, b, c\in X: a \cdot (a\rto c)\leq c$ and $(c\lto b)\cdot b \leq c$\ .
\eit
\end{defn}
For $1$ the unit of the monoid and setting ${}^*a:=a\rto 1$ and  $b^*:=1\lto b$, it then follows that  each Ajdukiewicz-Bar-Hillel pomonoid is a protogroup \cite{Lambek1}.
 
\begin{defn}\em
A \em residuated pomonoid \em \cite{Lambek0} is a pomonoid
\[
(X, \leq, (-\rto -), (-\lto-))
\]
such that for all $a, b\in X$ we have two Galois adjunctions:
\[
(a\cdot -)\dashv (a \rto -)\quad \mbox{\rm and}\quad  (-\cdot b)\dashv (- \lto b)\ ,
\]
that is, explicitly: 
\bit
\item $\forall a, b, c\in X: b\leq a \rto c \Leftrightarrow a\cdot b \leq c \Leftrightarrow a \leq c\lto b$\ ,
\eit
or, equivalently, using the alternative characterization of the adjunctions:
\bit
\item $\forall a, b, c\in X: a \cdot (a\rto c)\leq c\  \mbox{,} \ c\leq  a\rto (a \cdot c) \ \mbox{,} \   
(c\lto b)\cdot b \leq c \ \mbox{,} \  c\leq   (c \cdot b)\lto b$\ .
\eit
\end{defn}
From the second formulation in terms of four conditions it   immediately follows that each residuated pomonoid is a Ajdukiewicz-Bar-Hillel pomonoid.

But note in particular also that what we have here is a \em non-commutative \em generalization of eq.~(\ref{eq:Heyting}) which defined an implication connective, conjunction being replaced by the (evidently) non-commutative composition of words, and the implication $(-\Rightarrow-)$ now having a right-directed and right-directed counterpart, respectively $(-\rto -)$ and $(-\lto-)$.

\begin{defn}\em
A \em Grishin pomonoid \em \cite{Grishin} is a residuated pomonoid 
\[
(X, \leq, (-\rto -), (-\lto-), 0)
\]
with a special element $0\in X$ which is such that:
\bit
\item $0 \lto (a\rto 0)= a = (0\lto a) \rto 0$ \ .
\eit
\end{defn}  

With some work one can show that every Grishin pomonoid is a residuated pomonoid too e.g.~see \cite{Lambek1}.  Now, anticipating the following definition we can set: 
\[
{}^*a:=a\rto 0\ \ \ \ \mbox{ ,} \ \ \ \  a^*:=0\lto a \ \ \ \  \mbox{and} \ \ \ \  a+b:={}^*(b^*\cdot a^*)=({}^*b \cdot {}^*a)^*\ , 
\]
(where the last equality is quite easy to prove) and then we have:
\bit
\item $a\cdot {}^*a \leq 0$ \ , \ $1\leq {}^*a+ a$ \ , \ $b^*\cdot b \leq 0$ \ , \ $1\leq b+ b^*$ \ . 
\eit
We also have that $a\rto c = {}^*a + c$ and $c\lto b = c + b^*$.  So now the non-commutative implication resembles the Boolean implications that we discussed in Section \ref{sec:logicandproof}, the *-opertions playing the role of negation and the $+$-operation corresponding to the disjunction.

\begin{defn}\em 
A \em pregroup \em \cite{Lambek1} is a pomonoid   
\[
(X, \leq, (-)^*, {}^*(-)) 
\]
 where  $(-)^*: X\to X$ and ${}^*(-): X\to X$  are such that: 
\bit
\item $a \cdot {}^*a \leq 1\leq {}^*a\cdot a$ \ , \ $b^*\cdot b \leq 1\leq b \cdot  b^*$ \ .
\eit
\end{defn}

Hence, each  pregroup  is a Grishin pomonoid with $\cdot = +$ and $0=1$, and each pregroup is  a protogroup which satisfies two additional conditions.  In the case that $(-)^*={}^*(-)$ then we obtain a group with $a^*=a^{-1}$, hence the name `pre'-group.

Right, that was a bit of a zoo!  Still, there is a clear structural hierarchy.  Below the arrows represent the increase in equational content: 
\[
\xymatrix@=0.20in{
  \mbox{protogroup} \ar[d]  \\
\mbox{ABH pom.}\ar[d] \\
\mbox{residuated pom.}\ar[d] \\
\mbox{Grishin pom.}\ar[d] \\
 \mbox{pregroup} 
}
\]
There are four kinds of inequalities that play a key role, which either reduce or introduce types, and do this either in terms of (left/right) unary or a (left/right) binary connective.  Those in terms of a unary connective imply the corresponding ones involving a binary connectives. The following  table depicts these rules, with $e\in\{ 0,1\}$ and $\circ\in \{\cdot, +\}$.
\begin{center}
\begin{tabular}{c||c|c|}
 & unary connective  & binary connective \\
 \hline \hline
 type reduction
 & $a \cdot {}^*a \leq e$ \quad $b^* \cdot b \leq e$ 
 & $a \cdot (a\rto c)\leq c$ \quad $(c\lto b)\cdot b \leq c$\\
 \hline
 type introduction
 & $1\leq a\circ a^*$ \quad $1\leq  {}^* b\circ b$ 
 & $c\leq  a\rto (a \circ c)$ \quad $c\leq   (c \circ b)\lto b$\\
 \hline
\end{tabular}
\end{center}
 
\section{Processes witnessing existence}\label{sec:SMC}

In Sec.~\ref{sec:things} we observed that for the very general setting of things/systems and processes thereon orders witness the existence of a process between two things/systems.  In Sec.~\ref{sec: composition} we saw how  composition of things interacts with ordering.  On the other hand, in  Sec.~\ref{sec: composition} we also saw that  monoid structures  naturally arise when we consider process composition in the sense of one process happening \em after \em another process.

Here we will make the passage from order witnessing existence of processes to explicitly describing these processes. Since processes themselves also come naturally  with sequential composition, and if the systems on which these act also compose, we will obtain a structure with two interacting modes of composition. A dual perspective is that starting from a process structure on a fixed system, we allow for variation of the system:

\bigskip

\hrule

\[
\xymatrix@=0.80in{
\hspace{-1cm}\mbox{composable processes}\hspace{-1cm}\ar[rd]|{\mbox{allow varying systems}\hspace{1cm}} & & \hspace{-1cm}\mbox{ordered composable systems}\hspace{-1cm}\ar[ld]|{\hspace{1cm}\mbox{explicitly describe processes}}  \\
&\hspace{-3cm}\mbox{dually composable processes between composable systems}\hspace{-3cm}  &
}
\]

\hrule

\bigskip

Historically, this structure traces back to the work of Benabou \cite{Benabou} and MacLane \cite{MacLaneCoherence}. Symbolically, this is what one is dealing with:

\begin{defn}\label{def:SMC}\em 
A \em strict symmetric monoidal category \em ${\cal S}$ consists of:
\bit
\item a collection (typically a `class')  of things/systems $|{\cal S}|$\ ,
\item with a monoid structure $(S, \otimes, 1)$ thereon,
\eit
and for each pair $S, S'\in |{\cal S}|$ 
\bit
\item a collection (typically a `set')  of processes ${\cal S}(S, S')$\ , 
\eit
with two unital associative composition structures:
\bit 
\item $\forall S, S', S''\in |{\cal S}|$, $(-\circ -):{\cal S}(S, S')\times {\cal S}(S', S'')\to {\cal S}(S, S'')$\ , 
\item $\forall S, S', S'', S'''\in |{\cal S}|$, $(-\circ -):{\cal S}(S, S')\times {\cal S}(S'', S''')\to {\cal S}(S\otimes S'', S'\otimes S''')$\ .
\eit
We denote ${\cal S}(S, S')$ also as $f:S\to S'$. Explicitly, associativity and unitality are:
\bit
\item $\forall f: S\to S', g:S'\to S'', h:S''\to S'''$ we have $(h\circ g)\circ f=h\circ (g\circ f)$\ ,
\item $\forall f: S\to S', g:S''\to S''', h:S''''\to S'''''$ we have $(f\otimes g)\otimes h = f\otimes (g\otimes h)$\ ,
\item $\forall S\in  |{\cal S}|$ there exists an \em identity process \em  $1_S:S\to S$ which is such that for all $f: S'\to S, g:S \to S'$ we have that $1_S\circ f=f$ and $g\circ 1_S=g$\ ,
\item there exists an \em identity system \em ${\rm I}\in |{\cal S}|$ which is such that for all $S\in |{\cal S}|$ we have that  ${\rm I}\otimes S=S\otimes {\rm I}= S$.
\eit
These composition  structures moreover interact \em bifunctorialy\em, that is:
\bit
\item $\forall f: S\to S', f': S'\to S'', g: S'''\to S'''', g': S''''\to S'''''$ we have that:
\beq\label{eq:bifunct1}
(f'\circ f)\otimes(g'\circ g)= (f'\otimes g')\circ (f\otimes g)\ ,
\eeq
\item $\forall S, S'\in |{\cal S}|$ we have that: 
\beq\label{eq:bifunct2}
1_S\otimes 1_{S'}= 1_{S\otimes S'}\ .
\eeq
\eit
Finally,  we assume \em symmetry  \em, that is,  
\bit
\item $\forall S, S'\in |{\cal S}|$ there exists a \em symmetry  process \em $\sigma_{S, S'}:S\otimes S'\to S'\otimes S$, and  these are such that for all $f: S\to S', g:S''\to S'''$ we have that:
\[
\sigma_{S', S'''}\circ(f\otimes g)= (g\otimes f) \circ \sigma_{S, S''} \ .
\]
\eit
\end{defn}

What a mess, or better, what a syntactic mess!  The problem here is indeed of a syntactic nature.  The concept behind a strict symmetric monoidal category is intuitively obvious but one wouldn't get that intuition easily when reading the above definition. In fact, by presenting \em strict \em symmetric monoidal categories rather than general symmetric monoidal categories we already enormously  simplified the presentation.  In the strict case we assume associativity and unitality of the $\otimes$-connective on-the-nose, while, as discussed in the introduction, set-theory based mathematical models would typically be non-strict.  However, the physical reality itself is strict, which points at an inadequacy of its typical mathematical models. Let us recall  this physical conception of process theories.  

There is a notion of system, to which we from now on will refer as \em type\em, and for each pair of types of systems there are processes which take a system of the first type to the system of the second type.  These processes can be composed in two manners.  

The first manner is \em sequentially\em, that is one process taking place \em after \em  another process, the second process having the output type of the first process as its input type.  

The second manner is \em in parallel\em, that is, one process takes place \em while \em the other one takes place, without any constraints on the input and output types.  

Examples of particular systems and processes respectively are `nothing' (cf.~${\rm I}$ in Defn.~\ref{def:SMC} above), and `doing nothing' (cf.~$1_S$ in Defn.~\ref{def:SMC} above).  

But there is more, ... Even in the strict definition above much of the structure is about `undoing' unavoidable syntactic features.   

For example, symmetry simply means that there is no significance to the list-ordering of systems when writing $S\otimes S'$, that is, $S\otimes S'$ and $S'\otimes S$ describe one and the same thing, and we can use $\sigma_{S, S''}$ to pass from one description to the other.  

Of course, sometimes the order does matter, like in the case of words making up a sentence.  Swapping words evidently changes the meaning of a sentence, and in most cases would make it even meaningless as swapping words of different grammatical types would typically destroy the grammatical structure.  So here, rather than a strict symmetric monoidal category, we would consider a \em strict monoidal category\em, which boils down to Defn.~\ref{def:SMC} without the $\sigma_{S, S'}$-processes.

Turning our attention again on the `syntactic mess', even more striking than the role of symmetry in undoing the ordering of one-dimensional linear syntax, is the role played by eq.~(\ref{eq:bifunct1}).  To expose its `undoing'-nature we will need to change language, from one-dimensional linear syntax to two-dimensional pictures.  This will also bring us much closer to our desire of basing our conception of foundational mathematical structure on the idea of high-levelness in use, given that the pictorial presentation gets rid of the artifacts of set-theoretical representation, as illustrated in the introduction on the example of `three things'.  The diagrammatic language  indeed on-the-nose captures the idea that strict symmetric monoidal categories aim to capture, but then  still within a syntactic realm.   The study of these diagrammatic languages is becoming more and more prominent in a variety of areas of mathematics, including modern algebra and topology.

In the  two-dimensional pictures processes will be represented by boxes and the input- and output-systems by wires:
\begin{center}
\quad\qquad\epsfig{figure=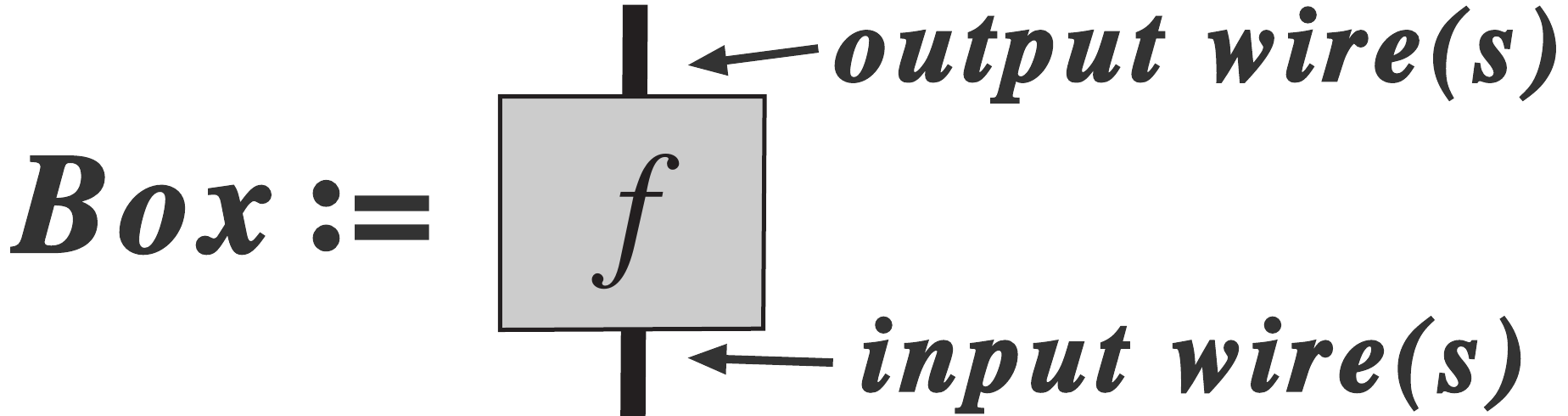,width=120pt}
\end{center}
We can then immediately vary systems by varying the number of wires:
\begin{center}
\ one system\    \quad $n$ \em sub\,\em-systems  \  \ \ \ \  \em no \em system \ \
\end{center}
\[
\underbrace{\epsfig{figure=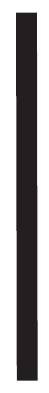,width=4pt}}_1\quad\ \ \ \qquad\underbrace{\epsfig{figure=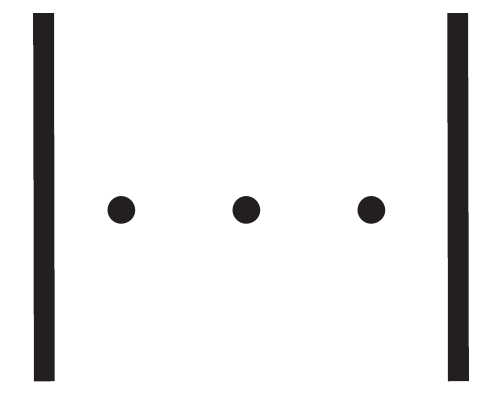,width=40pt}}_n\qquad\ \ \quad\underbrace{\epsfig{figure=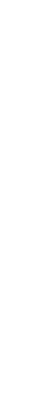,width=4pt}}_0
\]
The two compositions boil down to either connecting the output wire of one process to the input wire of the other, or by simply putting the processes side-by-side:
\[
\raisebox{-1.02cm}{\epsfig{figure=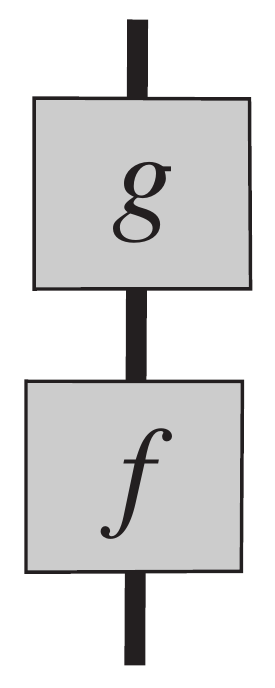,width=24pt}}
\qquad\qquad\qquad\qquad
\raisebox{-0.52cm}{\epsfig{figure=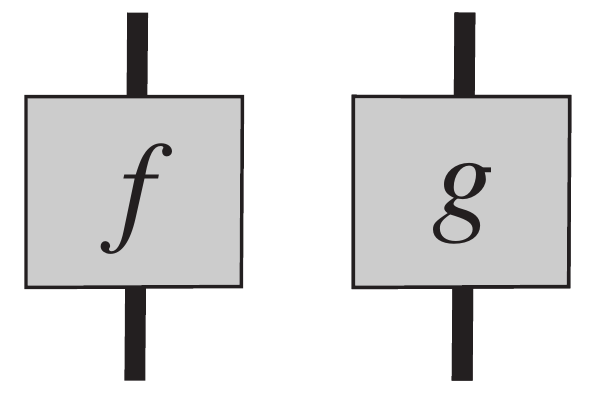,width=52pt}}
\]
Doing nothing is represented by a wire,  and nothing, evidently, by nothing.  

The rules of the game are: `only topology matters', that is, if two pictures are topologically equivalent, then they represent the same situation, e.g.:
\begin{center}
\epsfig{figure=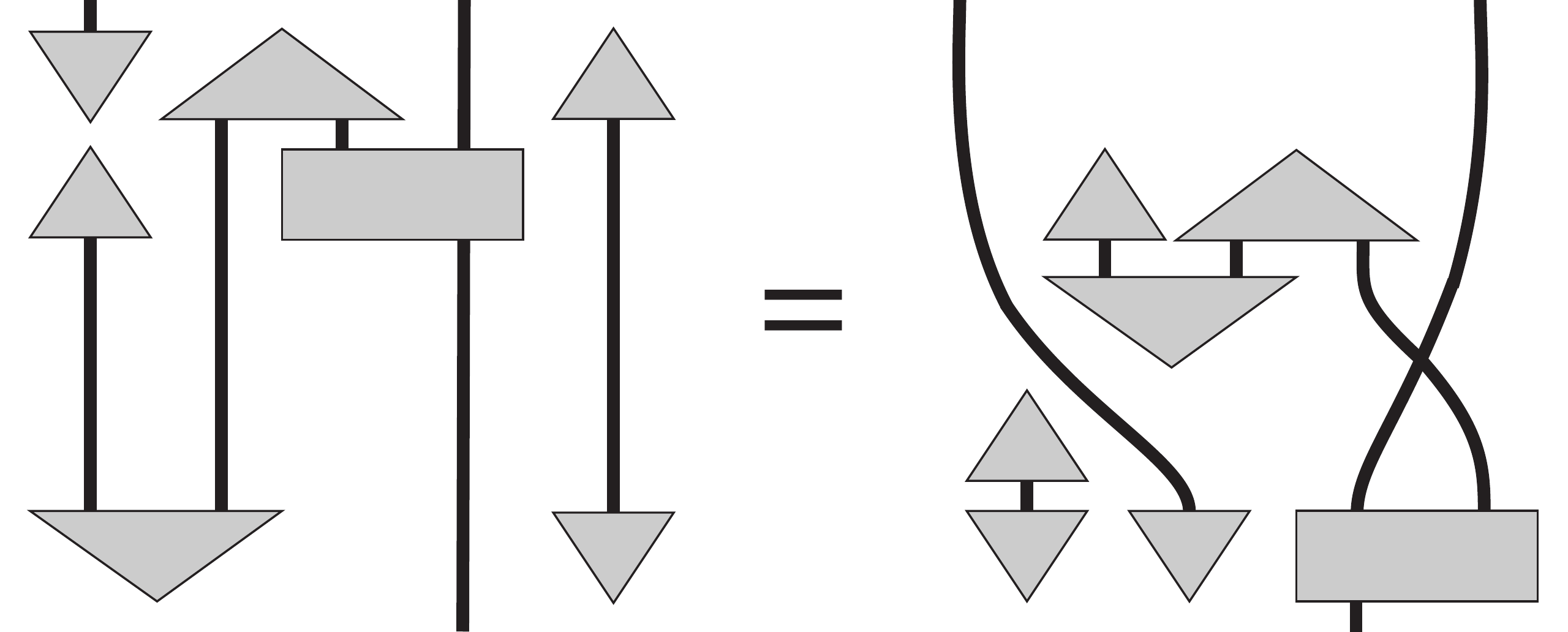,width=200pt}
\end{center}
Now, the nontrivial symbolic equation $(f'\circ f)\otimes(g'\circ g)= (f'\otimes g')\circ (f\otimes g)$ becomes:
\begin{center}
\epsfig{figure=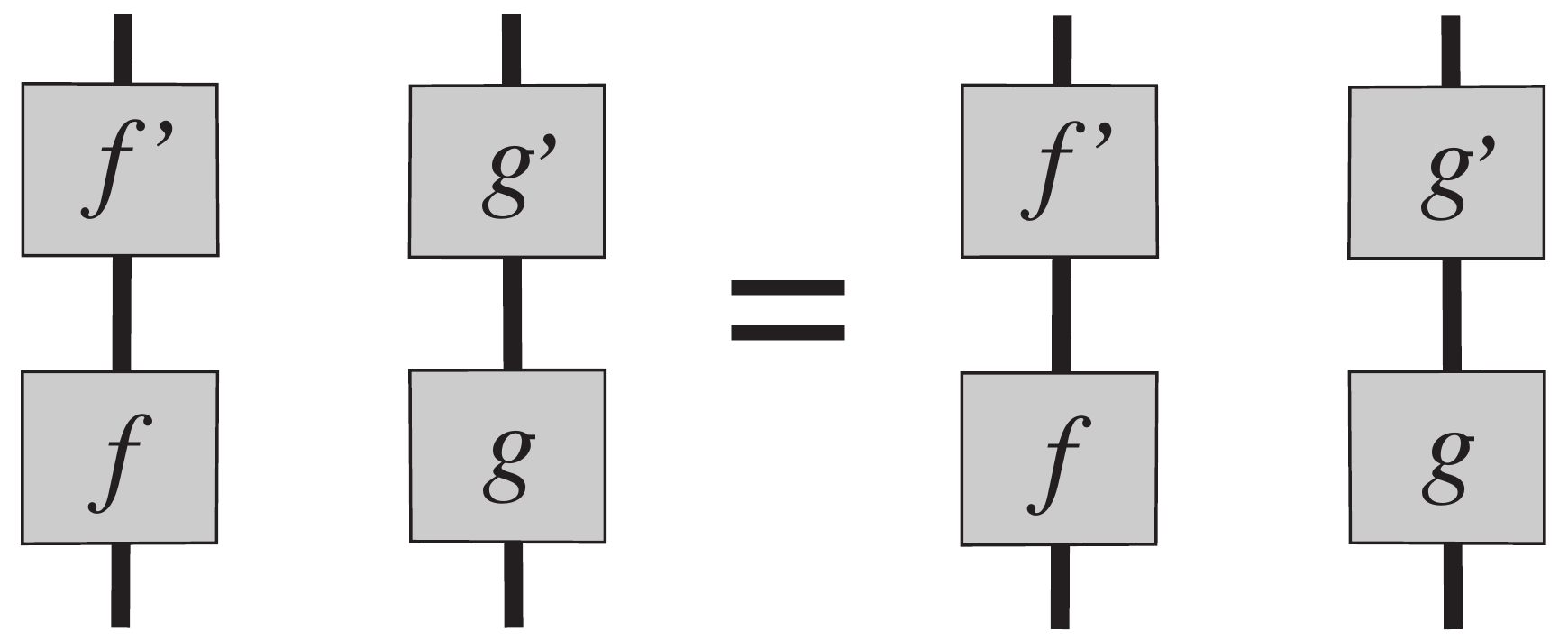,width=150pt}
\end{center}
i.e.~a tautology!  In other words, in this diagrammatic languages which more closely captures the mathematical structure of processes than its symbolic counterpart, essential symbolic requirements become vacuous. 

The reason is simple: there are two modes of composition that are in a sense `orthogonal', but one tries to encode them in a single dimension.  As a result, one needs to use brackets to keep the formulas well-formed, but these brackets obviously have no counterpart in reality.  This is where they would be in the pictorial language:
\begin{center}
\epsfig{figure=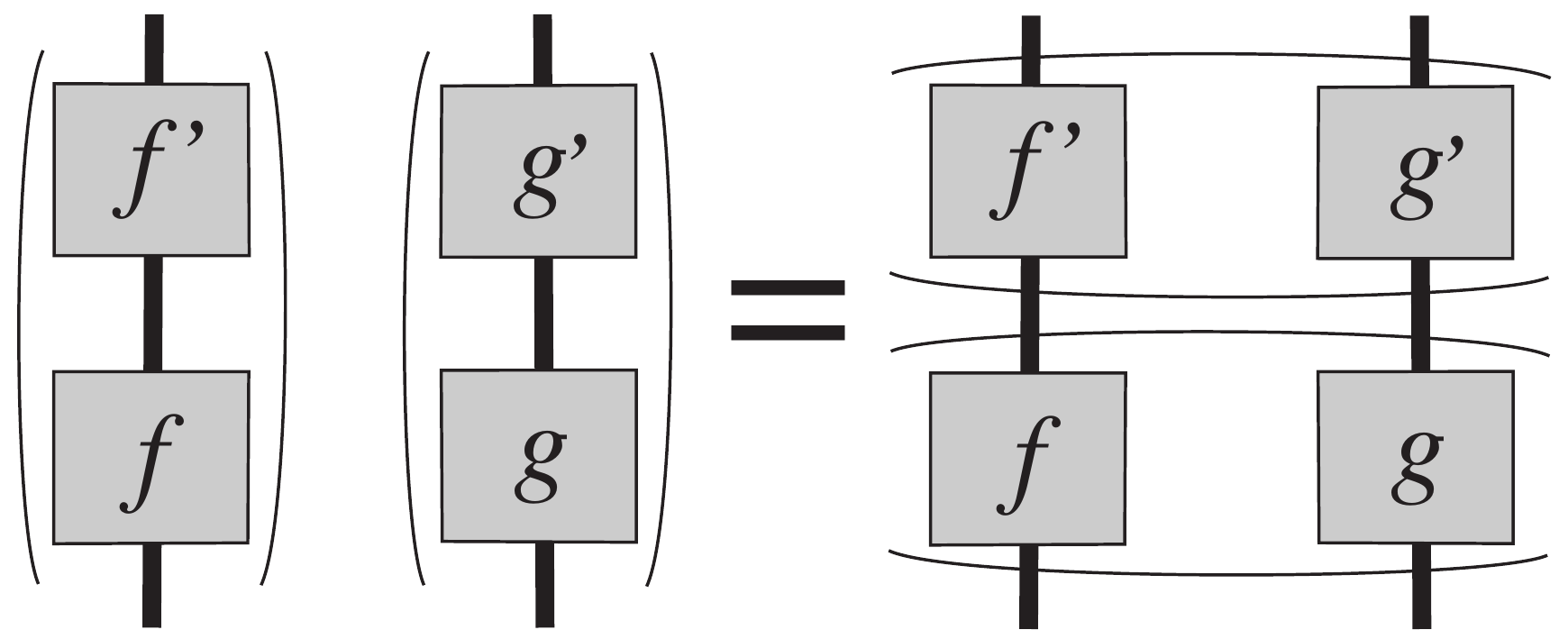,width=150pt}
\end{center}

A more detailed discussion of the upshot of graphical languages is in \cite{CatsII}.  The development of graphical languages for a variety of structures is an active area of research.  For a survey of the state of the art we refer to \cite{SelingerSurvey}.  Also useful in this context are \cite{KockBook, StreetBook, BaezLNP}.

\subsection{From-word-meaning-to-sentence-meaning processes}

We can blow up the pomonoids of Sec.~\ref{sec:catgram}  to full-blown proper (typically non-symmetric)  monoidal categories by replacing each relationship $a\leq b$ by  a collection ${\cal S}(a, b)$ of processes.  A residuated pomonoid then becomes a so-called \em bi-closed monoidal category\em, a Grishin pomonoid then becomes a category of which the symmetric counterpart is called a \em $*$-autonomous category \em \cite{Barr}, and a pregroup then becomes a category of which the symmetric counterpart is called a \em compact (closed) category \em \cite{Kelly, KellyLaplaza}.   

The non-symmetric case of  compact (closed) categories has been referred to as planar autonomous categories \cite{JS2, SelingerSurvey}, and the non-symmetric case of  $*$-autonomous categories as linearly distributive categories with negation \cite{Barr2, CockettSeely}.

These symmetric categories have themselves been studied in great detail, since closed symmetric monoidal categories capture \em multiplicative intuitionistic linear logic\em, while $*$-autonomous categories capture classical linear logic \cite{Girard, Seely}.  Compact closed categories model a degenerate logic which has found applications in quantum information processing \cite{AC1, RossThesis}.  We will discuss this a bit more in the next section. 

The following table summarizes the blow-up and symmetric restriction of some of the pomonoids  the arose when describing grammatical structure:

\begin{center}
\begin{tabular}{|c|c|c|} 
{\bf grammatical types} & {\bf monoidal category} & {\bf symmetric case} \\
 \hline 
residuated pomonoid & biclosed monoidal & closed symmetric monoidal\\
Grishin pomonoid & linearly distributive + negation & *-autonomous\\
pregroup & planar autonomous & compact closed \\
\hline
\end{tabular}
\end{center}

While, as explained in Sec.~\ref{sec:catgram},  the ordered structures capture how grammatical structure evolves when composing words, the categorical structures capture how the meaning of words transforms into the meaning of sentences \cite{CSC,GrefSadr}.  The diagrammatic representation of these categories then explicitly shows how meaning `flows' within sentences.   Here is an example of such a `meaning flow' taken from \cite{CSC}:
\begin{center}
\epsfig{figure=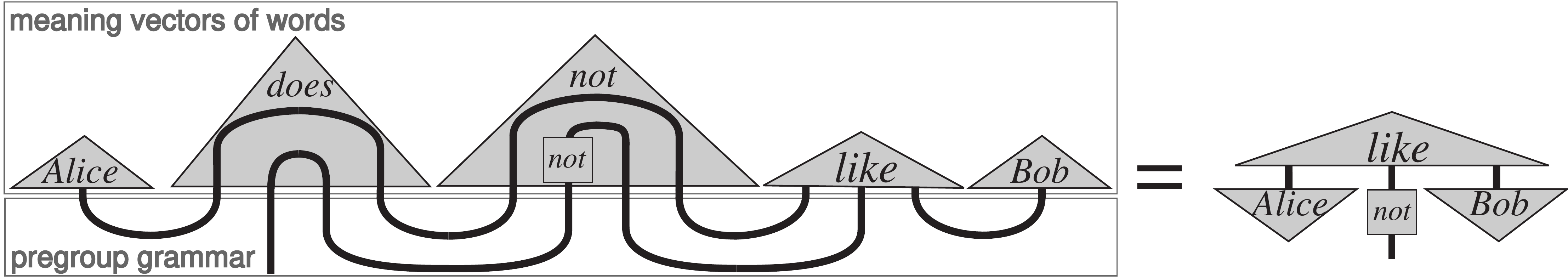,width=360pt}
\end{center}
The verb ``like'' receives the subject ``Alice'' via the flow of meaning through ``does'' and ``not'', and also receives the subject ``Bob''.  Then it produces the meaning of the sentence ``Alice does like Bob'', which then is negated by the not-box.

The reader may verify that this particular explicit representation of meaning flow requires at the grammatical level the full structure of a Grishin pomonoid. 


\subsection{Discipline transcending process structures}

In fact, very similar pictures when modeling  information flows in quantum protocols, in work that inspired the above one on language meaning \cite{AC1, ContPhys, QLog}:
\begin{center}
\epsfig{figure=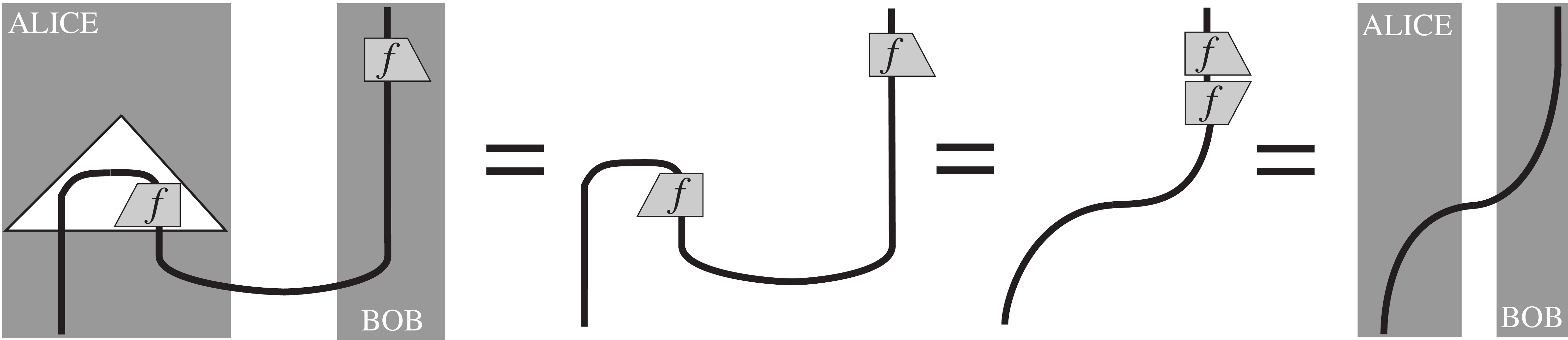,width=360pt}
\end{center}
Here Alice and Bob share a Bell-state (the `cup'-shaped wire), then Alice performs an operation that depends on a discrete variable $f$, and also Bob does an operation that depends on $f$.  The end-result is a perfect channel between Alice and Bob.  This protocol is known as quantum teleportation \cite{Tele}.

Another example is probabilistic Bayesian inference \cite{CSBayes}:
\begin{center}
\epsfig{figure=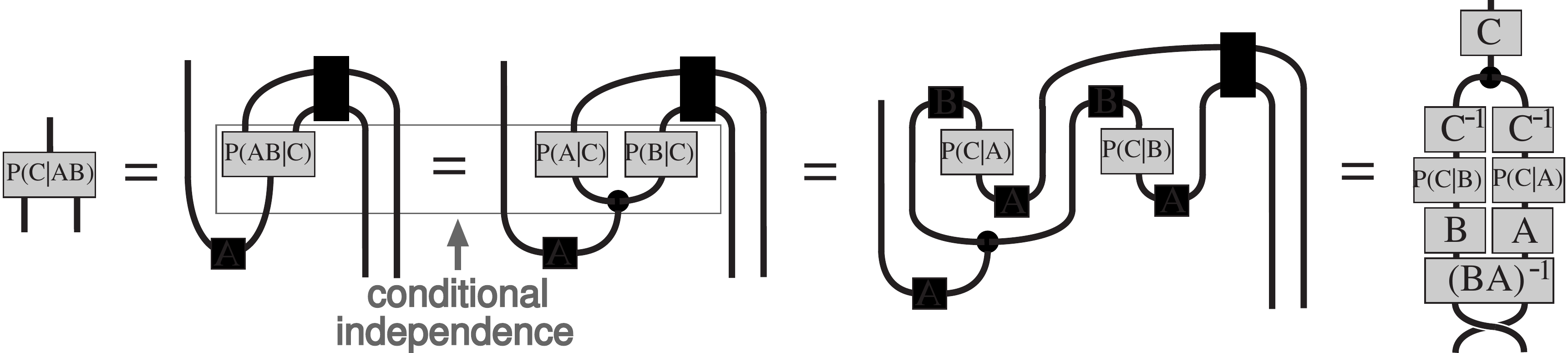,width=360pt}
\end{center}
Rather than just wires we now also have `dots' connecting several wires.  These structures are also highly relevant in the quantum mechanical applications (e.g.~see \cite{CD2}), and since recently, also in the linguistic applications, where they play the role of bases.  They are also key to quantum algebra \cite{Majid, StreetBook} and topological quantum field theory \cite{KockBook, Turaev}.  The precise connections between these uses have yet to been fully  explored.

%

\section{Closing}  

We started by discussing ordering on things, to which we then adjoined composition of things, and by passing from existence of processes that take certain things to other things, to explicitly representing them, we ended with a structure which is most naturally represented, not by a syntactic, but by a diagrammatic language.  

We gave examples of applications of these structures in a wide range of disciplines.  The message that we tried to pass to the reader is that these structure are very basic, in that they appear in such a huge range of situations when taking a high-level perspective.  

It would therefore be natural to, before doing anything fancy, to give these structures a privileged status. Then, as a first next step, other well-understood structure could come in as, for example, so-called categorical enrichment \cite{KellyEnriched, BorceuxStubbe}. Evidently, it would be nice to take this even further, and see how far one could ultimately get by setting up this style of structural hierarchy driven by high-levelness of actual phenomena.

\section{Acknowledgements}

The author is enjoys support for the British Engineering and Physical Sciences Research Council, from the Foundational Questions Institute, and from the John Templeton Foundation.  We thank Robin Cockett and Robert Seely for filling some holes in our background knowledge, and Chris Heunen's for pointing out typos.

\bibliographystyle{plain}

\bibliography{LinguisticsBookarXiv} 
\end{document}